\documentclass[graybox, vecphys]{svmult}

\usepackage{type1cm}        
\usepackage{makeidx}         
\usepackage{graphicx}        
\usepackage{multicol}        
\usepackage[bottom]{footmisc}

\usepackage{newtxtext}       %
\usepackage[varvw]{newtxmath}       

\makeindex             

\usepackage{cleveref}
\usepackage{bigfoot}
\usepackage{wasysym}

\begin{document}

\title*{Domain decomposition solvers for operators with fractional interface perturbations}
\author{Miroslav Kuchta}
\institute{Miroslav Kuchta \at Simula Research Laboratory, Kristian Augusts gate 23, 0164 Oslo \email{miroslav@simula.no}}
%
%
\maketitle

\abstract{
  Operators with fractional perturbations are crucial components for robust preconditioning
  of interface-coupled multiphysics systems. However, in case the perturbation is strong,
  standard approaches can fail to provide scalable approximation of the inverse, thus compromising
  efficiency of the entire multiphysics solver. In this
  work, we develop efficient and parameter-robust algorithms for interface-perturbed operators based
  on the non-overlapping domain decomposition method. As preconditioners for the
  resulting Schur complement problems we utilize (inverses of) weighted sums of fractional powers of the
  interfacial Laplacian. Realization of the preconditioner in terms of rational approximation is
  discussed. We demonstrate performance of the solvers by numerical examples including application
  to coupled Darcy-Stokes problem.
}

\section{Introduction}\label{sec:intro}
Mathematical models featuring interaction of physical systems across a common interface
describe numerous phenomena in engineering, environmental sciences and medicine. Here the
large variations in material coefficients or wide ranges of temporal/spatial scales at
which the phenomena can be studied demand parameter-robust solution algorithms. In \cite{boon2022robust, boon2022parameter}    
such algorithms were recently developed for Darcy-Stokes and Biot-Stokes models   
by establishing uniform stability of the respective problems in (non-standard) parameter-dependent
norms. In particular, the authors show that in order to attain robustness, mass conservation at the
interface $\Gamma$ of the porous domain $\Omega$ must be accounted for in the functional setting,
leading to control of the porous pressure $p$ in the norm $\lVert p \rVert_{\Omega}$ such that
\begin{equation}\label{eq:norm}
\lVert p \rVert_{\Omega}^2 = \lVert K^{1/2} \nabla p \rVert^2_{0, \Omega} + \lVert \mu^{-1/2} p \rVert^2_{-1/2, \Gamma}.
\end{equation}
Here $\lVert \cdot \rVert_{k, D}$ denotes the standard norm of Sobolev space $H^k(D)$ on domain $D$.
The coefficients $K, \mu>0$ are due to material properties, namely the permeability of
the porous medium and the fluid viscosity.

By operator preconditioning, the choice of norm \eqref{eq:norm} yields a Riesz map preconditioner
$b \mapsto x$ defined by solving the problem
\begin{equation}\label{eq:riesz}
   -K\Delta_{\Omega} x + \mu^{-1}(-\Delta_{\Gamma})^{-1/2} x = b.
\end{equation}
Note that the operator in \eqref{eq:riesz} contains a bulk part $-\Delta_{\Omega}$ and
an interface part $(-\Delta_{\Gamma})^{-1/2}$, which, from the point of view topological dimension
of the underlying domains, can be viewed as a lower order \emph{perturbation}.  

\begin{figure}[t]
  \centering
  \includegraphics[height=0.25\textwidth]{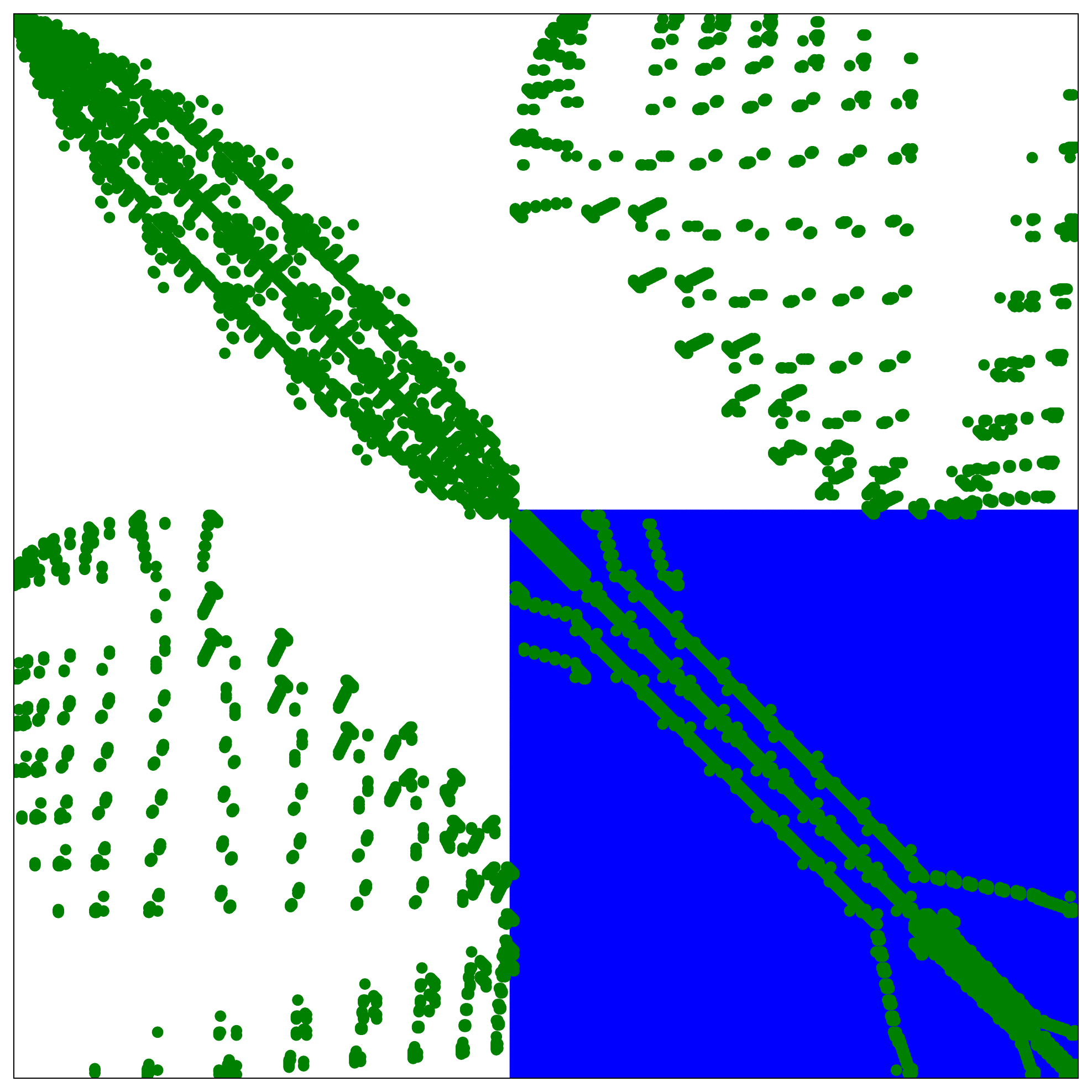}
  \hspace{30pt}
  \includegraphics[height=0.25\textwidth]{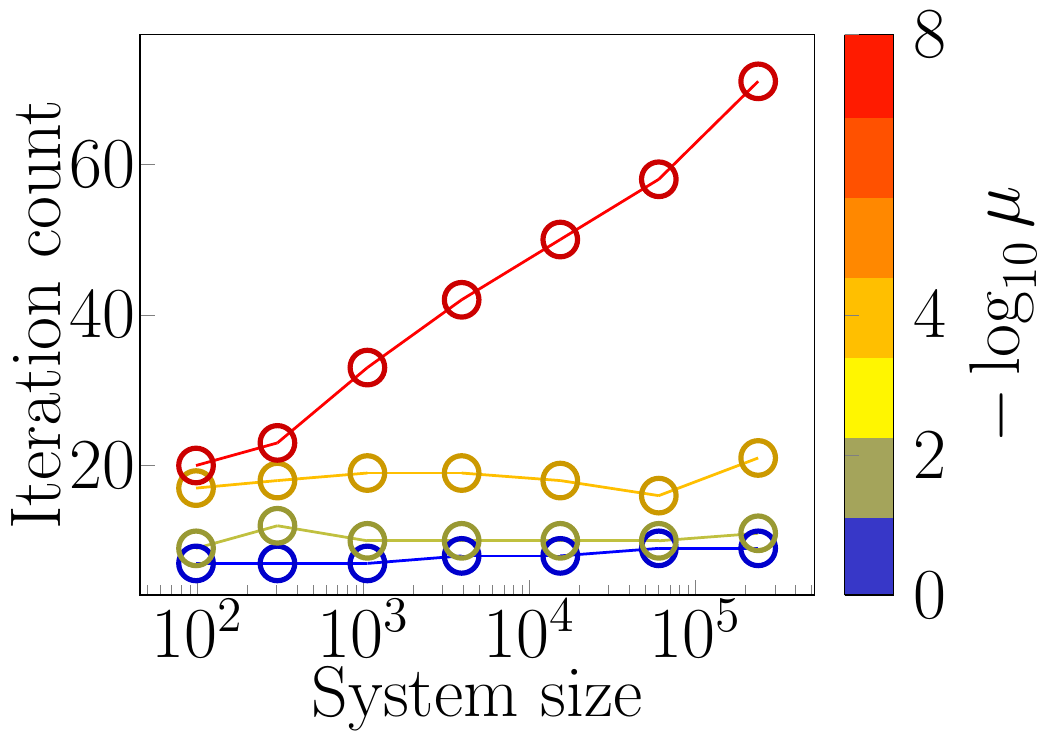}
  \vspace{-5pt}
  \caption{
    (Left) Sparsity pattern of the operator in \eqref{eq:riesz} on $\Omega=(0, 1)^3$ with $\Gamma\subset\partial\Omega$.
    Interface perturbation leads to dense block (in blue) which is challenging for sparse LU
    solvers. (Right) Number of PCG iterations under mesh refinement when solving \eqref{eq:riesz} on
    $\Omega=(0, 1)^2$ with $\Gamma\subset\partial\Omega$ and AMG \cite{hypre}   
    preconditioner. In both  case case $K=1$ and the problems are discretized by continuous linear Lagrange
    ($\mathbb{P}_1$) elements. 
  }
  \label{fig:sparsity_amg}
\end{figure}

Efficiency of the block-diagonal Darcy/Biot-Stokes preconditioners \cite{boon2022robust, boon2022parameter} hinges   
on performant solvers for \eqref{eq:riesz}. However, the problem might not be amenable
to standard (generic, black-box) approaches in particular in case when the fractional interface
perturbation becomes dominant. We illustrate this behaviour in \Cref{fig:sparsity_amg}
where \eqref{eq:riesz} is solved by preconditioned conjugate gradient (PCG) method with
algebraic multigrid (AMG) preconditioner. Indeed, the number of iterations increases with the
weight of the perturbation term and, worryingly, for large enough values mesh-independence
is lost. 

Non-overlapping domain decomposition (DD) is a solution methodology which has been successfully
applied to number of challenging problems including coupled multiphysics systems e.g.
\cite{wietse, discacciati2018optimized, vassilev2014domain, galvis2009feti}.  
A key component of the method are then the algorithms for the problems arising at the interface
which can be broadly divided into two categories. In FETI or BDDC variants (see e.g. \cite{arioli2013discrete} and references therein) the solvers utilize
suitable auxiliary problems on the \emph{subdomains}. To develop tailored solvers for operators with
fractional interface perturbation we here follow an alternative approach \cite{arioli2013discrete} and address the problem
directly at the \emph{interface}. In particular, we shall construct preconditioners for the resulting
Steklov-Poincar{\'e} operators using sums of fractional order interfacial operators which include contribution
due to the DD and the perturbation (which is only localized at the interface).


\section{Domain decomposition solvers}\label{sec:solver}
We shall consider solvers for \eqref{eq:riesz} in a more general setting. To this end, let 
$\Omega\subset\mathbb{R}^d$, $d=2, 3$ be a bounded domain with boundary $\partial\Omega$
and $\Gamma\subseteq\partial\Omega$. Moreover, let $V=V(\Omega)$, $Q=Q(\Gamma)$ be a pair of
Hilbert spaces with $V'$, $Q'$ being their respective duals and let $R:V\rightarrow Q'$ be
a restriction operator. For $b\in V'$ we are then interested in solving
\begin{equation}\label{eq:generic}
\mathcal{A}x = b\text{ in }V'\quad\text{with}\quad\mathcal{A}=A_{\Omega} + \gamma R' B^{-1}_{\Gamma} R,
\end{equation}
where $\gamma\geq 0$ and $A_{\Omega}:V\rightarrow V'$ is some symmetric operator coercive on $V$ while
$B_{\Gamma}:Q\rightarrow Q'$ is assumed to induce an inner product on $Q$.
Note that the norm operator in \eqref{eq:riesz} is a special case of \eqref{eq:generic} with
$V=H_{0}^{1}(\Omega)$, $Q=H^{1/2}(\Gamma)$, $R$ the trace operator
and $A_{\Omega}=-K\Delta_{\Omega}$ while $B_{\Gamma}=(-\Delta_{\Gamma})^{1/2}$.


To formulate our non-overlapping domain decomposition approach for \eqref{eq:generic}, we follow \cite{arioli2013discrete} and 
decompose $V=V_0\oplus V_{\Gamma}$ where $V_0=\left\{v \in V; Rv = 0\right\}$. Assuming that $V_{\Gamma}$ can be
identified with $Q$ we observe that the operator $\mathcal{A}$ takes a block structure
\begin{equation}\label{eq:block_op}
  \arraycolsep=2pt
  \def\arraystretch{1.2}  
  \mathcal{A} = \left(\begin{array}{c|c}
    A^{00}_{\Omega} & A^{0i}_{\Omega}\\
    \hline
    A^{i0}_{\Omega} & A^{ii}_{\Omega}
    \end{array}
    \right) +
    \gamma\left(\begin{array}{c|c}
    0 & 0\\
    \hline
    0 & \tilde{B}^{-1}_{\Gamma}
    \end{array}
    \right),
\end{equation}
which we exploit to design a preconditioner for $\mathcal{A}$. Specifically, 
under the assumption that $A^{00}_{\Omega}$ is invertible, let us define the DD preconditioner
\begin{equation}\label{eq:dd_preconditioner}
  \arraycolsep=2pt
  \def\arraystretch{1.2}  
  \mathcal{B} = \left(\begin{array}{c|c}
    I^{00}_{\Omega} & -(A^{00}_{\Omega})^{-1}A^{0i}_{\Omega}\\
    \hline
    0 & I^{ii}_{\Omega}
    \end{array}
  \right)
\left(\begin{array}{c|c}
    A^{00}_{\Omega} & 0\\
    \hline
    0 & S_{\Gamma}
    \end{array}\right)^{-1}
\left(\begin{array}{c|c}
    I^{00}_{\Omega} & 0\\
    \hline
    -A^{i0}_{\Omega}(A^{00}_{\Omega})^{-1} & I^{ii}_{\Omega}
    \end{array}\right).
\end{equation}
Here $I^{00}_{\Omega}:V_0\rightarrow V_0$ and $I^{ii}_{\Omega}:V_{\Gamma}\rightarrow V_{\Gamma}$
are identity operators on the respective subspaces while $S_{\Gamma}$ is spectrally equivalent
to the DD Schur complement/Steklov-Poincar{\'e} operator. We note that preconditioner
\eqref{eq:dd_preconditioner} preserves symmetry of the original problem \eqref{eq:generic} as we
target PCG solvers. However, with Krylov methods which do not require symmetry a more efficient
triangular variant of the preconditioner is sufficient.

Our main contribution is an observation that for problems with interface perturbations, the Schur
complement preconditioner $S_{\Gamma}$ in \eqref{eq:dd_preconditioner} takes the form
\begin{equation}\label{eq:schur_precond}
  S_{\Gamma} = L^s_{A} + \gamma L^t_{B},
\end{equation}
for some constants $s, t\in \mathbb{R}$ and symmetric, positive-definite operators $L_A$, $L_B$
depending on the regularity of $\mathcal{A}$, the restriction operator and the perturbation.
In particular, the structure of the preconditioner reflects
the two contributions to the Schur complement; the decomposition $V=V_0\oplus V_{\Gamma}$ applied to
operator $A_{\Omega}$ yields $L^s_{A}$ while $L^t_{B}$ is due to the perturbation.

Motivated by the initial example \eqref{eq:riesz} we shall in the following focus on problems for
which $-1 < s, t < 1$ and $L_A$, $L_B$ are spectrally equivalent to $L=-\Delta_{\Gamma}+I_{\Gamma}$. However,
we highlight that the operators might in general differ by their boundary conditions (which for $L_A$ reflect
boundary conditions on $\partial\Omega\setminus\Gamma$).


Assuming that $(A^{00}_{\Omega})^{-1}$ can be efficiently computed, the main challenge for scalability of
preconditioner \eqref{eq:dd_preconditioner} is an efficient realization of (an approximate) inverse of \eqref{eq:schur_precond}.
Upon discretization, the operators $L_{A}^s$, $L_{B}^t$ can be approximated by eigenvalue factorization\footnote{
For $L:Q\rightarrow Q'$ let $L_h$ be the matrix realization of the operator in the basis of some finite dimensional
approximation space $Q_h$, $n=\text{dim}\,Q_h$. Moreover, let $M_h$ be the mass matrix, i.e. matrix realization of the inner product
of the Lebesgue space $L^2$ on $Q_h$. Assuming $L$ is symmetric, positive definite, the factorization
$L_hU_h=M_hU_h\Lambda_h$, $U^T_hM_hU_h=\text{Id}$ holds where $\Lambda_h$ is a diagonal matrix of eigenvalues while the corresponding
eigenvectors constitute the columns of matrix $U_h$. We then define
\begin{equation}\label{eq:gevp}
  L^s_h = (M_h U_h) \Lambda^s_h (M_h U_h)^{T}.
\end{equation}
Note that for $L=(-\Delta + I)$ and $f\in Q$ represented in $Q_h$ by interpolant with coefficient vector
$f_h\in\mathbb{R}^n$ the function $f_h\mapsto f_h\cdot L^s_hf_h $ represents an approximation
of the square of the Sobolev norm $\lVert f \rVert^2_{s}$. 
}.
However, this approach suffers from cubic scaling. For the specific case of $L^{1/2}$ a more efficient strategy with
improved scaling is applied in \cite{arioli2013discrete} based on the Lanczos process while, more recently, \cite{harizanov2022rational} prove
that rational approximations (RA) lead to non-overlapping DD methods with linear scaling. Building on this observation to
obtain order optimal solvers for the perturbed problem \eqref{eq:generic} we follow \cite{budisa2022rational} where rational
approximations\footnote{Referring to the definitions in \eqref{eq:gevp} the RA construct approximate solutions
$x\in Q$ satisfying $\alpha L^sx + \beta L^tx=b$, $b\in Q'$ in the finite dimensional space $Q_h$ via a solution operator
\begin{equation}\label{eq:RA}
c_0M^{-1}_h + \sum^m_{k=1}c_i(L_h + p_k M_h)^{-1}.
\end{equation}
Here, $c_i$ and $p_i\geq 0$ are respectively the residues and the poles of the rational approximation $f_{\text{RA}}$
to function $f: x\rightarrow (\alpha x^s + \beta x^t)^{-1}$. Importantly, the number of poles $m$ does not depend on the
dimensionality of $Q_h$ and is instead determined by the accuracy $\epsilon_{\text{RA}}$ of the RA, i.e.
$\lVert f - f_{\text{RA}}\rVert\leq \epsilon_{\text{RA}}$. We refer to \cite{harizanov2022rational, budisa2022rational} and references
therein for more details. 
}
were developed for Riesz maps with norms induced by sum operator
$\alpha L^s + \beta L^t$ with $\alpha, \beta\geq 0$. In particular, this setting fits our Schur complement operator
\eqref{eq:schur_precond} if constant material properties \emph{and} suitable boundary conditions are prescribed
on $\mathcal{A}$ in \eqref{eq:generic}.
%

\section{Model problem}\label{sec:numerics_base}
We shall illustrate performance of the domain decomposition preconditioner \eqref{eq:dd_preconditioner}
using a model interface-perturbed problem: Find $x\in V=H^1(\Omega)$ such that
\begin{equation}\label{eq:model_problem}
K(-\Delta_{\Omega}+I_{\Omega})x + \gamma(-\Delta_{\Gamma}+I_{\Gamma})^tx = b\text{ in }V',
\end{equation}
where $K>0$, $\gamma\geq 0$ and $-1<t<1$. Here $\Omega=(0, 1)^d$, $d=2, 3$ and $\Gamma=\partial\Omega$.
We note that this choice maximizes the size of the interface. At the same times, it enables the RA-favorable
setting of $L_A=L$, $L_B=L$, $L=-\Delta_{\Gamma}+I_{\Gamma}$ in \eqref{eq:schur_precond}. 
Following \cite{arioli2013discrete} the DD Schur complement of the operator $A_\Omega=K(-\Delta_{\Omega}+I_{\Omega})$
in \eqref{eq:generic} is spectrally equivalent to fractional operator $KL^{1/2}$. In turn we apply
preconditioner \eqref{eq:dd_preconditioner} with $S_{\Gamma}=KL^{1/2} + \gamma L^t$. 

In the numerical experiments we consider $H^1$-conforming finite element spaces 
$V_h\subset V$ constructed in terms of $\mathbb{P}_1$ elements. Consequently, the matrix realization of the
fractional interface perturbation then reads $\gamma T^{T}_h L^t_h T_h$ where matrix $L^t_h$ is defined
in \eqref{eq:gevp} and $T_h$ is a discrete trace operator such that $T_h\phi=\sum^n_{j=1}l_j(\phi|_{\Gamma})\psi_j$
for any $\phi\in V_h$ and $\phi_j$, $l_j$, $j=1, \cdots, n$ being respectively the basis functions and
degrees of freedom (point evaluations) of the discrete trace space $V_{\Gamma, h}=Q_h$ built likewise using $\mathbb{P}_1$ elements.

The linear systems due to discretized \eqref{eq:model_problem} shall be solved by PCG solver using
our DD preconditioner \eqref{eq:dd_preconditioner} which now requires inverse of the linear system due to
$K L_h^{1/2} + \gamma L_h^t$. Here we shall either apply the eigenvalue realization \eqref{eq:gevp} (which allows
for closed form evaluation of the exact inverse) or the approximate inverse due to RA, see \eqref{eq:RA}.
To put focus on the Schur complement the leading block $(A^{00}_{\Omega})^{-1}$ in \eqref{eq:dd_preconditioner} will
be computed exactly by LU factorization. For results with approximate inverse of $A^{00}_{\Omega}$ we refer to \Cref{rmrk:lazy}.
Finally, the PCG solver is always started from 0 initial vector and terminates upon reducing the preconditioned
residual norm by factor $10^{10}$.

\begin{figure}[t]
  \centering
  \includegraphics[width=\textwidth]{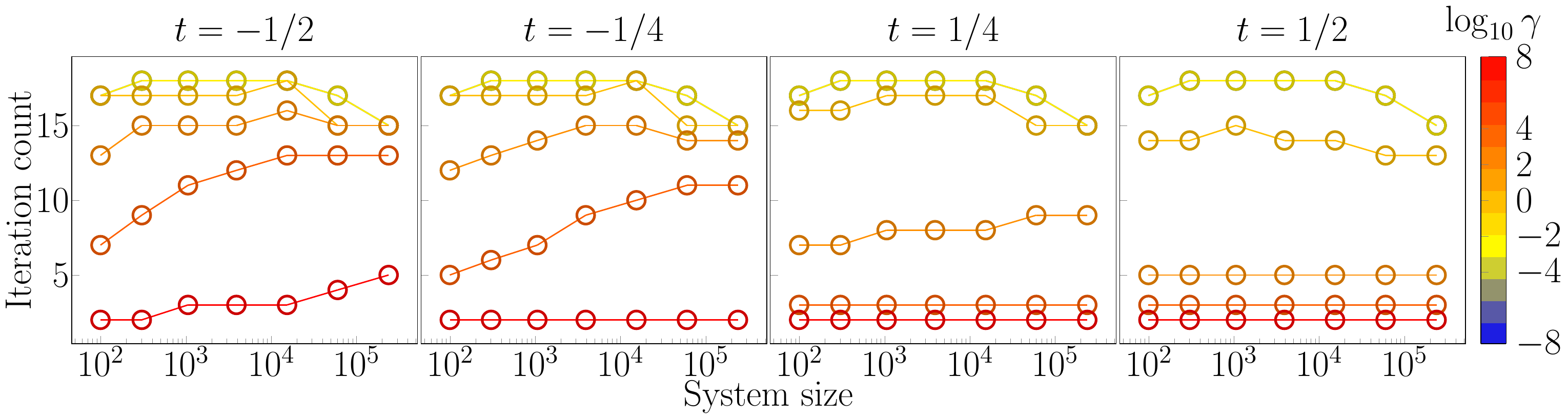}
  \vspace{-19pt}      
  \caption{
    PCG iterations when solving \eqref{eq:model_problem} on $\Omega=(0, 1)^2$ and preconditioner
    \eqref{eq:dd_preconditioner} with $S_{\Gamma}=KL^{1/2} + \gamma L^t$. Problem is discretized by $\mathbb{P}_1$
    elements. Blocks of the preconditioner are here computed exactly.
  }
  \label{fig:base_operator_2d}
\end{figure}

\begin{figure}
  \centering
  \includegraphics[width=\textwidth]{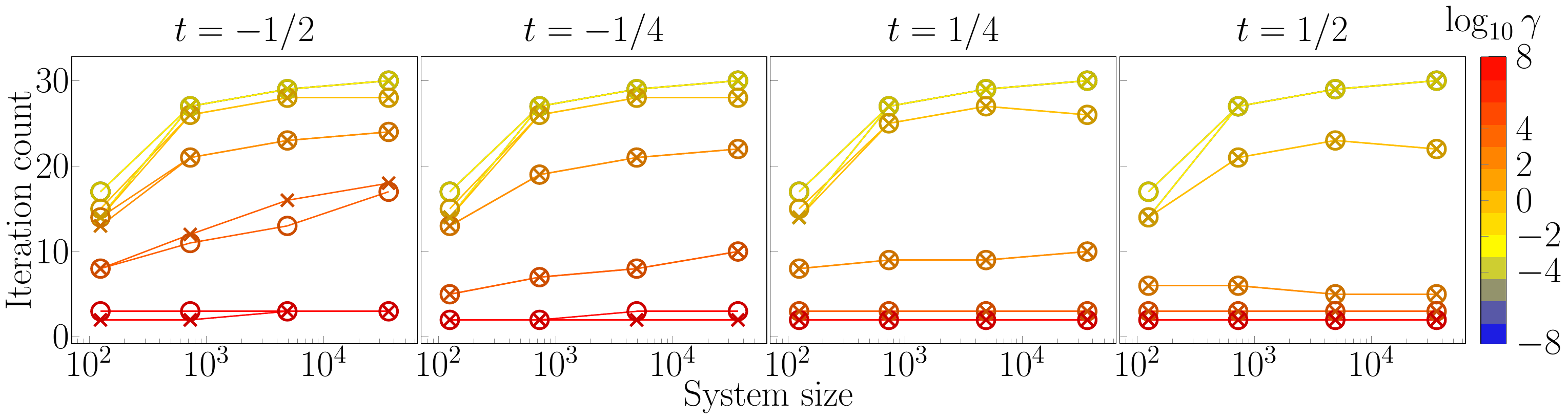}
  \vspace{-19pt}
  \caption{
    PCG iterations when solving \eqref{eq:model_problem} on $\Omega=(0, 1)^3$ and preconditioner
    \eqref{eq:dd_preconditioner} with $S_{\Gamma}=KL^{1/2} + \gamma L^t$. Problem is discretized by $\mathbb{P}_1$
    elements. Leading block of preconditioner is computed exactly. 
    Results with realization of the Schur complement preconditioner by RA with tolerance $\epsilon_{\text{RA}}=10^{-14}$
    are depicted by ($\ocircle$) markers while ($\times$) markers correspond to definition
    via the eigenvalue problem \eqref{eq:gevp}.
  }
  \label{fig:base_operator_3d}
\end{figure}

We summarize performance of the DD preconditioner in \Cref{fig:base_operator_2d} and \Cref{fig:base_operator_3d} which
consider \eqref{eq:model_problem} with $\Omega=(0, 1)^2$ and $\Omega=(0, 1)^3$ respectively. It can be seen that the
PCG convergence is in general bounded in mesh size, fractionality $t$ and the perturbation strength $\Gamma$.
Important for the scalability of \eqref{eq:dd_preconditioner} is the fact that iteration counts with RA realization of the
Schur complement preconditioner practically match the exact inverse of $S_{\Gamma}$. We remark that the chosen tolerance
of $\epsilon_{\text{RA}}=10^{-14}$ yields roughly $m=20$ poles in \eqref{eq:RA}. The computation setup in $3d$ then leads to
linear systems with $<6200$ unknowns at the interface. 


\begin{figure}
  \centering
  \includegraphics[height=0.25\textwidth]{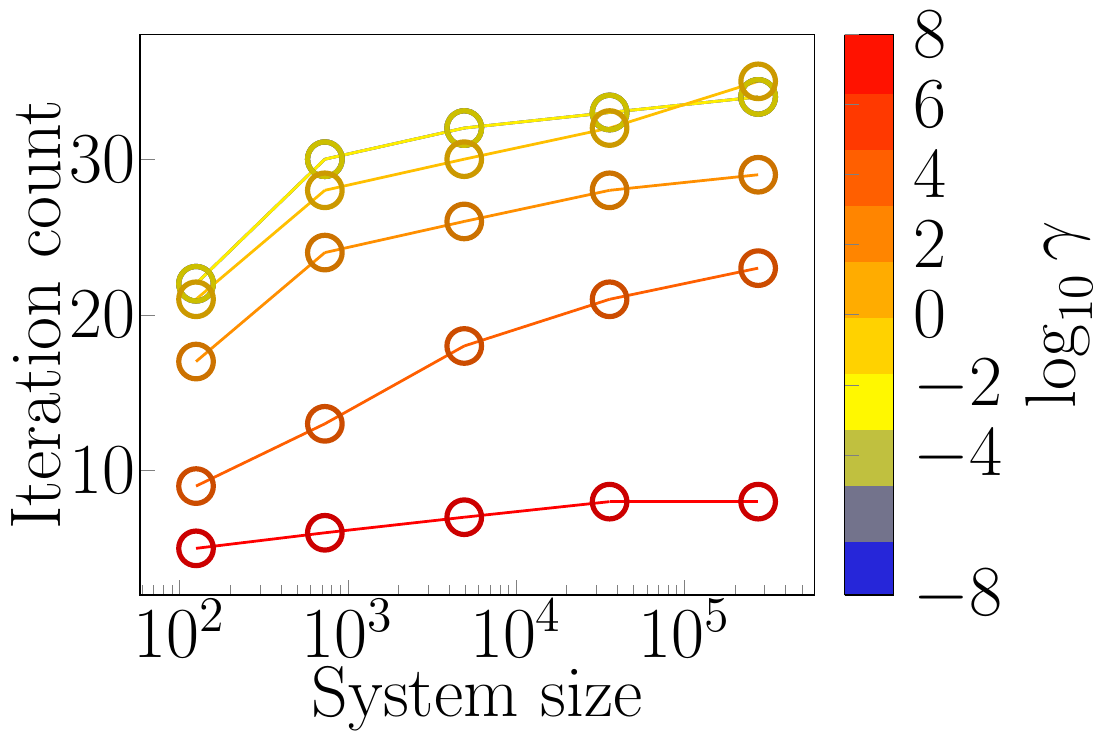}
  \hspace{30pt}
  \includegraphics[height=0.25\textwidth]{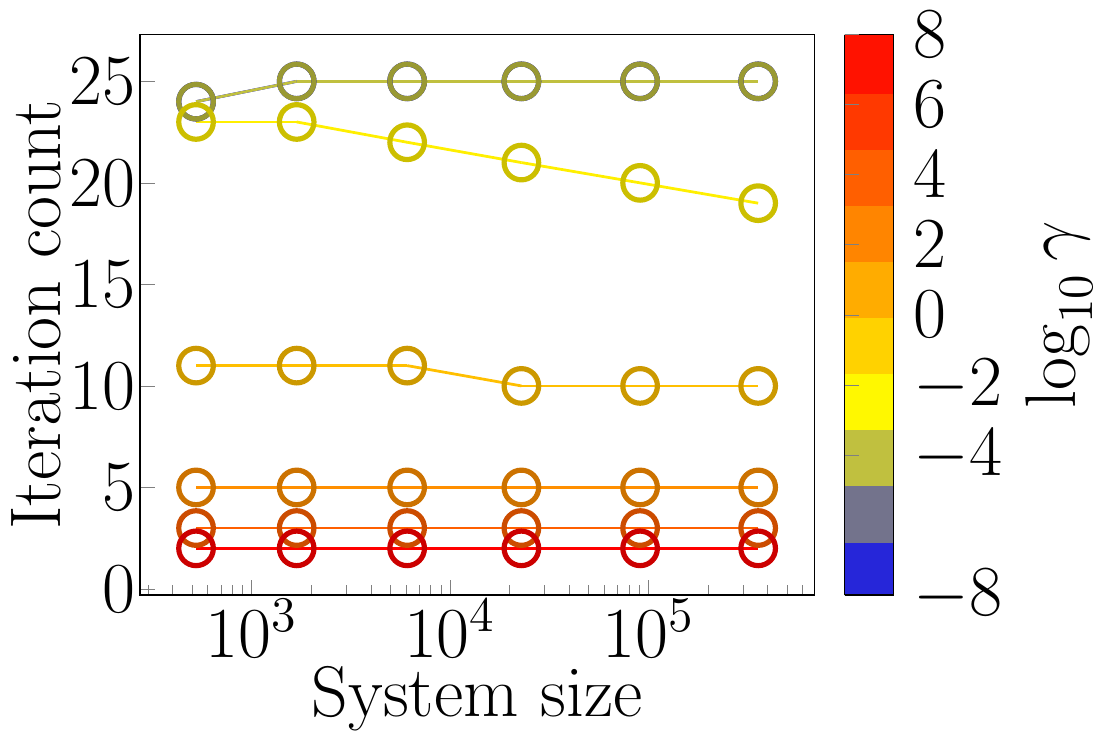}
  \vspace{-9pt}        
  \caption{
    PCG iterations for computing the inverse of fractional perturbed operators by
    preconditioner \eqref{eq:dd_preconditioner}. (Left) The operator is \eqref{eq:model_problem}
    with $t=-1/2$ and $\Omega=(0, 1)^3$. Both the operator and preconditioner are
    evaluated by RA. On the finest refinement level $\text{dim}\,V_{\Gamma, h}=24\cdot 10^3$. (Right) 
    Operator \eqref{eq:graddiv} is considered with $S_{\Gamma}=KL^{-1/2}+\gamma I_{\Gamma}$ in the
    Schur complement \eqref{eq:schur_precond}.
  }
  \label{fig:graddiv_and_lazy}
\end{figure}

\begin{remark}[Evaluation of the operator in \eqref{eq:model_problem}]\label{rmrk:lazy}
  In numerical experiments shown in \Cref{fig:base_operator_2d}, \Cref{fig:base_operator_3d} the
  operator $\mathcal{A}$ in \eqref{eq:model_problem} utilized the eigenvalue decomposition \eqref{eq:gevp}
  of $L^t_h$. This realization restricts the size of $\Gamma$ or $\text{dim}\,Q_h$
  that can be considered. However, action of the perturbation can instead be computed via RA leading
  to evaluation of $\mathcal{A}$ with optimal complexity and enabling large scale problems. In \Cref{fig:graddiv_and_lazy}
  we revisit \eqref{eq:model_problem} with $t=-1/2$, $\Omega=(0, 1)^3$ and RA used both for the operator and
  the preconditioner \eqref{eq:dd_preconditioner}. Moreover, the leading block of the preconditioner is approximated
  by a single V-cycle of AMG \cite{hypre}.
  The number of PCG iterations then appears to be
  bounded in the mesh size and the parameter $\gamma$. As before, $\mathbb{P}_1$ elements were used for discretization.
  %
\end{remark}

\begin{remark}[Application to \vec{H}(\text{div})-elliptic problem]
  Preconditioners \eqref{eq:dd_preconditioner} are not limited to $H^1$-elliptic problems.
  To illustrate this fact we consider $V=\vec{H}(\text{div}, \Omega)$, $\Omega=(0, 1)^2$ and a
  variational problem induced by bilinear form due to operator $\mathcal{A}$
  \begin{equation}\label{eq:graddiv}
    \langle \mathcal{A} \vec{u}, \vec{v} \rangle = \int_{\Omega}K(\vec{u}\cdot\vec{v} +
    \nabla\cdot\vec{u}\nabla\cdot\vec{v}) + \gamma\int_{\partial\Omega}\vec{u}\cdot{\vec{\nu}}\vec{v}\cdot{\vec{\nu}}
    \quad\forall \vec{u},\vec{v}\in V.
  \end{equation}
  We observe that $\mathcal{A}$ falls under the template problem \eqref{eq:generic}. In order to
  apply the domain-decomposition preconditioner we then require a preconditioner for the DD Schur complement
  due to $(A^{00}_{\Omega})^{-1}$ where $A^{00}_{\Omega}$ is here the operator $K(I - \nabla\nabla\cdot)$
  on $\vec{H}_0(\text{div}, \Omega)$. Motivated by \cite{babuvska2003mixed}, we shall to this end consider the operator
  $L_A=KL^{-1/2}$ such that $S_{\Gamma}$ in \eqref{eq:dd_preconditioner} is defined as $S_{\Gamma}=KL^{-1/2}+\gamma I_{\Gamma}$.
  For numerical experiments the system is discretized by lowest order Brezzi-Douglas-Marini elements which lead
  to the discrete trace space $V_{\Gamma, h}=Q_h$ of discontinuous piecewise-linear functions on trace mesh $\Gamma_h$.
  Robustness of the resulting preconditioner is shown in \Cref{fig:graddiv_and_lazy}.
\end{remark}

\section{Darcy-Stokes preconditioning}\label{sec:numerics_ds}
\begin{figure}
  \centering
  \includegraphics[height=0.25\textwidth]{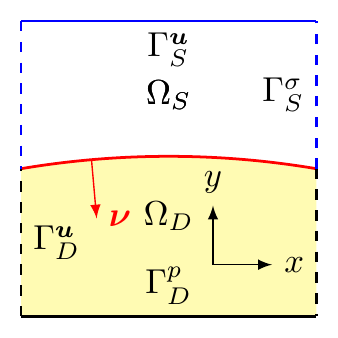}
  \hspace{5pt}
  \includegraphics[height=0.25\textwidth]{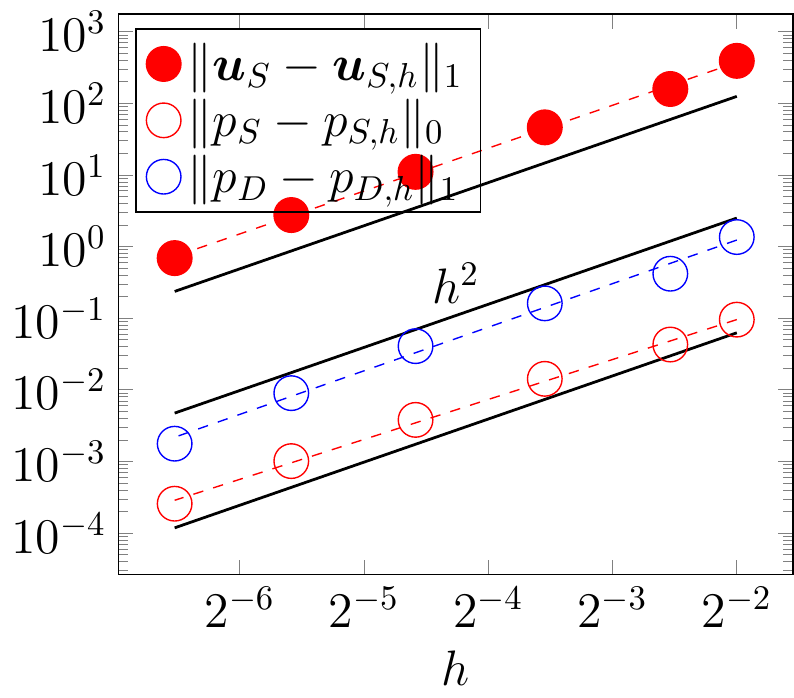}
  \hspace{5pt}  
  \includegraphics[height=0.25\textwidth]{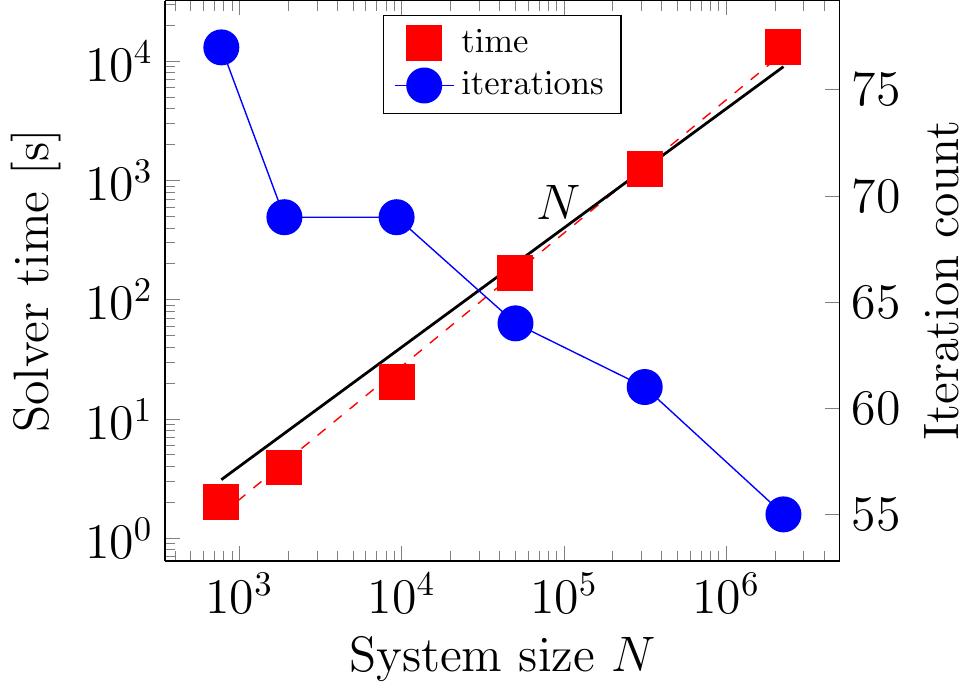}
  \vspace{-10pt}
  \caption{
    (Left) Computation domain is obtained by extrusion of the pictured geometry. The interface $\Gamma$, being part of a circle arc,
    is curved. No-slip and traction conditions are prescribed
    on $\Gamma^{\vec{u}}_S$ and $\Gamma^{\sigma}_S$ respectively. Darcy pressure is prescribed on $\Gamma^{p}_D$.
    (Center) Error convergence study performed using the $3d$ setup. With $\mathbf{\mathbb{P}}_2$-$\mathbb{P}_1$-$\mathbb{P}_2$
    elements, optimal quadratic rates are observed in all the variables.
    (Right) The solver time (including preconditioner setup and FGMRes runtime) scales linearly with the problem size.
  }
  \label{fig:darcy_stokes_3d}
\end{figure}

We finally apply the proposed non-overlapping DD solvers to realize preconditioners for
the coupled Darcy-Stokes problem with Darcy problem in the primal form \cite{discacciati2018optimized}. That is,
assuming bounded domains $\Omega_S$, $\Omega_D\subset\mathbb{R}^d$, $d=2, 3$ sharing a common interface $\Gamma$
(cf. \Cref{fig:darcy_stokes_3d}) we seek Stokes velocity and pressure $\vec{u}_S$, $p_S$ and Darcy pressure $p_D$ such that
\begin{equation}\label{eq:darcy_stokes}
  \begin{aligned}
  -\nabla\cdot\sigma(\vec{u}_S, p_S) = \vec{f}_S \text{ and } \nabla\cdot \vec{u}_S &= 0  &\text{ in }\Omega_S, \\
  -\nabla\cdot K\nabla p_D &= f_D  &\text{ in }\Omega_D, \\
  \vec{u}_S\cdot \vec{\nu} + K\nabla p_D\cdot \vec{\nu} &= 0  &\text{ on }\Gamma, \\
  -\vec{\nu} \cdot\sigma(\vec{u}_S, p_s)\cdot\vec{\nu}  - p_D &= 0   &\text{ on }\Gamma, \\
  -P_{\vec{\nu} }\left(\sigma(\vec{u}_S, p_S)\cdot\vec{\nu}_S\right) - \alpha\mu K^{-1/2}P_{\vec{\nu}} \vec{u}_S &= 0   &\text{ on }\Gamma,
  \end{aligned}
\end{equation}
where $P_{\vec{\nu}}$ is the tangential trace operator $P_{\vec{\nu}}\vec{u}=\vec{u}-(\vec{u}\cdot\vec{\nu})\vec{\nu}$
and $\sigma(\vec{u}, p)=\mu\Delta\vec{u}-p\text{Id}$. In addition
to the previously introduced coefficients $K$, $\mu>0$ the model also includes the Beavers-Joseph-Saffman parameter $\alpha\geq 0$.
We close the system by prescribing suitable boundary conditions to be discussed shortly. 

We consider \eqref{eq:darcy_stokes} with a parameter-robust block diagonal preconditioner \cite{boon2022robust}
\begin{equation}\label{eq:darcy_stokes_precond}
  \mathcal{B} = \text{diag}\left(-\mu\Delta + \alpha \mu K^{-1/2} P'_{\vec{\nu}} P_{\vec{\nu}},
  \mu^{-1}I,
  -K\Delta + \mu^{-1}(-\Delta_{\Gamma})^{-1/2}
  \right)^{-1}.
\end{equation}
Observe that both the first and the final block in \eqref{eq:darcy_stokes_precond} are of the form of the
interface-perturbed operators \eqref{eq:generic}. However, for simplicity we shall here set $\alpha=0$
and only focus on the pressure preconditioner. In particular, to efficiently approximate \eqref{eq:riesz}
we shall perform few PCG iterations with the DD preconditioner \eqref{eq:dd_preconditioner} using
$S_{\Gamma}=KL^{1/2}+\mu^{-1}L^{-1/2}$. We note that the interface operator is thus identical to the one utilized
in robust preconditioning of mixed Darcy-Stokes model \cite{holter2020robust}.

To illustrate performance of the preconditioner \eqref{eq:darcy_stokes_precond} we consider \eqref{eq:darcy_stokes}
in a $3d$ domain pictured in \Cref{fig:darcy_stokes_3d} and set\footnote{
  Due to computational demands we did not perform parameter-robustness study for $d=3$.
  In the subsequent Appendix we demonstrate robustness using a $2d$ version of the geometry in \Cref{fig:darcy_stokes_3d}.
}
$K=10^{-2}$, $\mu=10^{-4}$. Using discretization by
$\mathbf{\mathbb{P}}_2$-$\mathbb{P}_1$-$\mathbb{P}_2$ elements the linear system is solved by
preconditioned Flexible GMRes (FGMRes) starting from 0 initial guess and relative tolerance $10^{-10}$ on the preconditioned
residual norm. The Darcy-Stokes preconditioner is
then realized by applying single AMG V-cycle for the Stokes blocks while the Riesz map of the Darcy pressure \eqref{eq:riesz}
is approximated by PCG solver using \eqref{eq:dd_preconditioner} and running with a relative tolerance of $10^{-4}$.
The DD preconditioner uses RA with tolerance $\epsilon_{\text{RA}}=10^{-14}$ and AMG for the leading block of \eqref{eq:dd_preconditioner}.
With this setup the scalability study summarized in \Cref{fig:darcy_stokes_3d} reveals that the proposed solver is order optimal.
We note that the computations are run in serial.


%

\vspace{-20pt}
\section*{Appendix}\label{sec:appendix}
\addcontentsline{toc}{section}{Appendix}
We verify parameter-robustness of Darcy-Stokes preconditioner \eqref{eq:darcy_stokes_precond}
which utilizes \eqref{eq:dd_preconditioner}-preconditioned CG solver for the Darcy pressure block
using the ($2d$) geometry and boundary condition setup depicted in \Cref{fig:darcy_stokes_3d}. Here, a relative
tolerance of $10^{-5}$ is prescribed in PCG, while (i) either the RA in \eqref{eq:dd_preconditioner} is applied
with $\epsilon_{\text{RA}}=10^{-12}$ and AMG for the leading block or (ii) the respective problems are solved with
\eqref{eq:gevp} and LU. For simplicity and to put focus on the presented DD approach the Stokes velocity and
pressure blocks of \eqref{eq:darcy_stokes_precond} are computed exactly with LU. FGMRes settings are identical
to \Cref{sec:numerics_ds}.

\begin{figure}
  \centering
  \includegraphics[width=\textwidth]{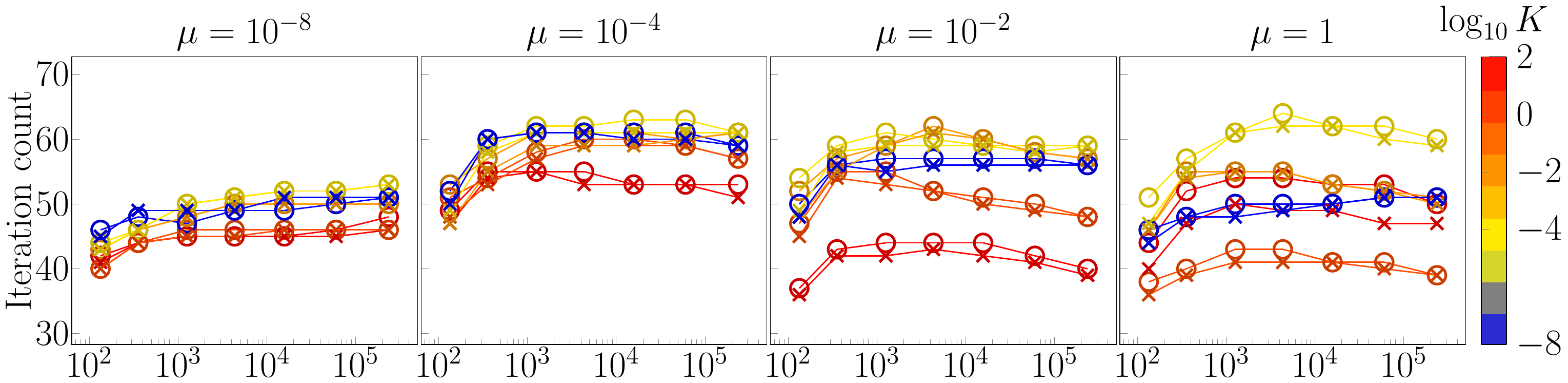}
  \includegraphics[width=\textwidth]{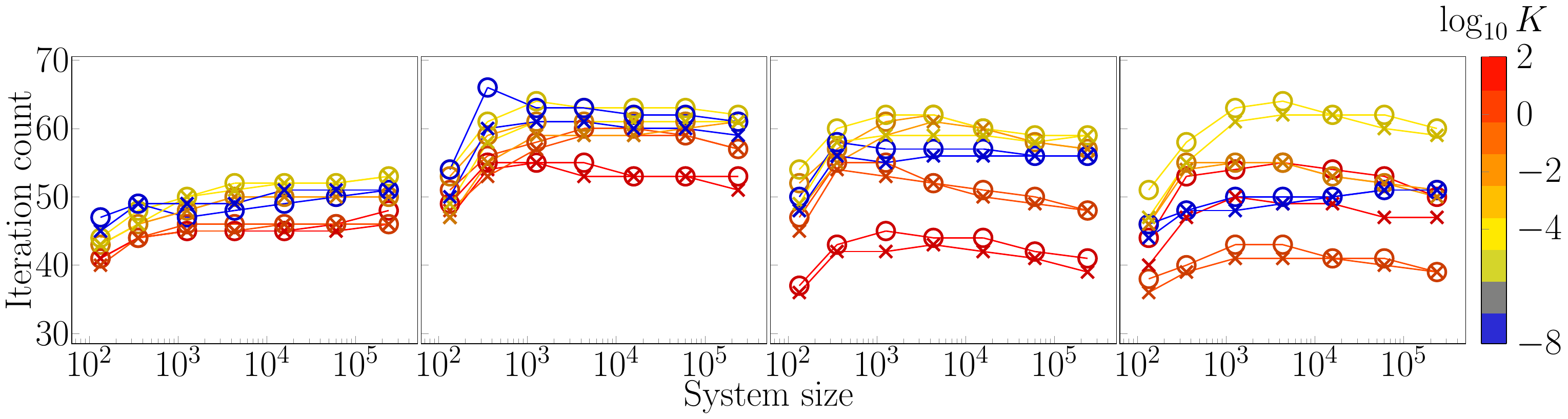}
  \vspace{-20pt}
  \caption{
    FGMRes iterations for Darcy-Stokes problem \eqref{eq:darcy_stokes} on $2d$ geometry
    in \Cref{fig:darcy_stokes_3d}. Beavers-Joseph-Saffman parameter is set to $\alpha=3$. Discretization
    by $\mathbf{\mathbb{P}}_2$-$\mathbb{P}_1$-$\mathbb{P}_2$ elements. Darcy-Stokes preconditioner \eqref{eq:darcy_stokes_precond}
    uses LU for Stokes velocity and pressure blocks. 
    (Top) Two different realizations of \eqref{eq:riesz} are compared.
    Results depicted in ($\times$) markers are obtained with \eqref{eq:riesz} computed by LU.
    With ($\ocircle$) markers the Riesz map \eqref{eq:riesz} is approximated
    by PCG using relative tolerance $10^{-5}$ and DD preconditioner \eqref{eq:dd_preconditioner}.
    Within \eqref{eq:dd_preconditioner} $(A^{00}_{\Omega})^{-1}$ is computed by LU and the Schur complement uses \eqref{eq:gevp}. 
    (Bottom) DD preconditioner \eqref{eq:dd_preconditioner} within PCG solver for \eqref{eq:riesz}
    using LU and \eqref{eq:gevp} for the blocks ($\times$) is compared with a scalable version
    utilizing AMG, respetively RA ($\epsilon_{\text{RA}}=10^{-12}$) to approximate the blocks of \eqref{eq:dd_preconditioner}.
  }
  \label{fig:darcy_stokes_2d}
\end{figure}

\begin{figure}
  \centering
  \includegraphics[width=\textwidth]{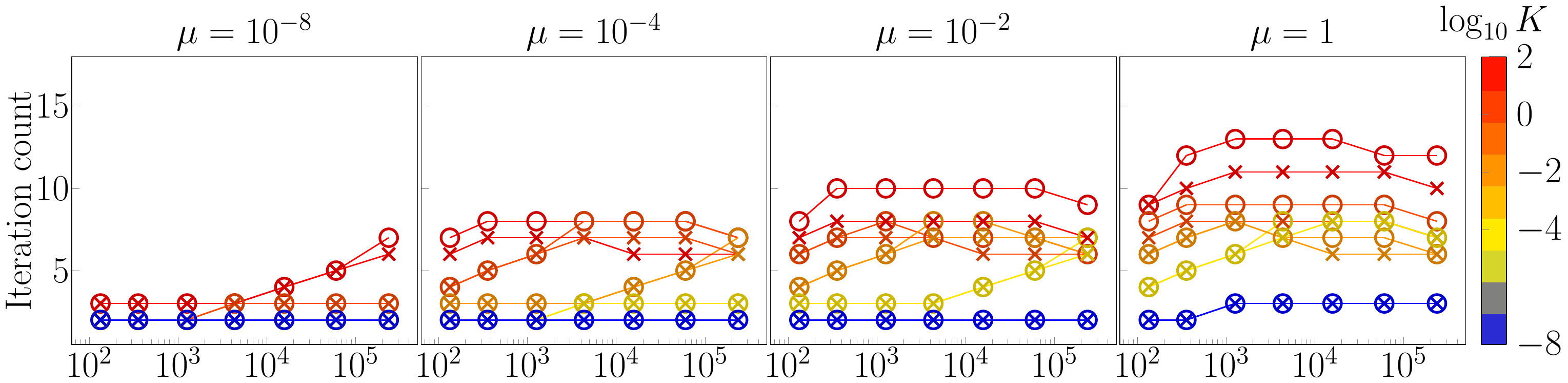}
  \includegraphics[width=\textwidth]{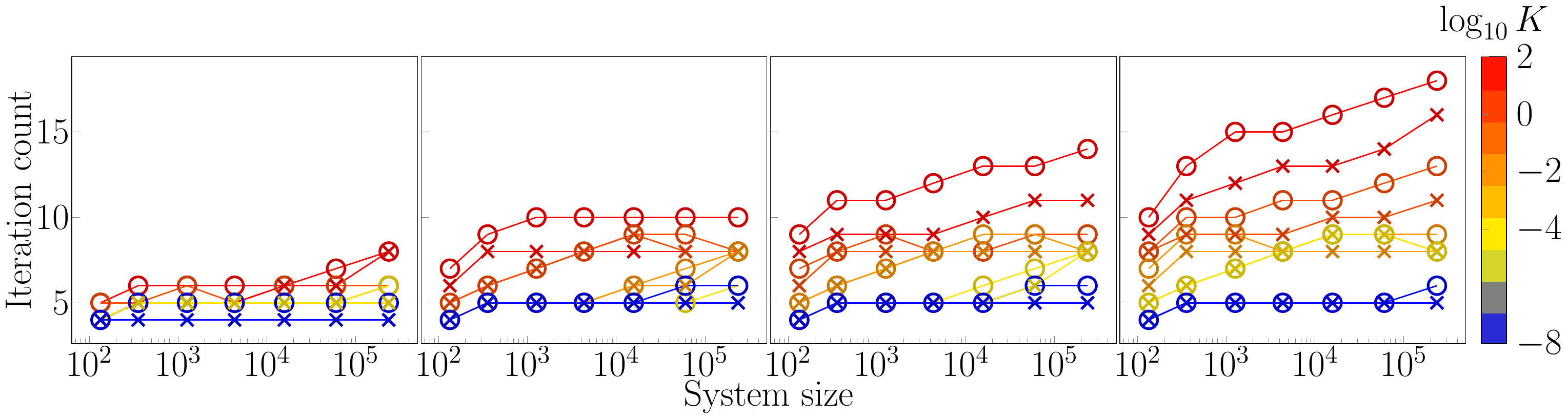}
  \vspace{-20pt}  
  \caption{
    PCG convergence inside application of Darcy-Stokes preconditioner in \Cref{fig:darcy_stokes_2d}.
    (Top) DD preconditioner uses LU  to realize $(A^{00}_{\Omega})^{-1}$ and \eqref{eq:gevp} for $S^{-1}_{\Gamma}$.
    (Bottom) The blocks of \eqref{eq:dd_preconditioner} are approximated by AMG V-cycle and RA ($\epsilon_{\text{RA}}=10^{-12}$)
    respectively. In both plots, 
    values in ($\ocircle$) markers represent the maximum number of PCG iterations needed for convergence
    when applying the outer preconditioner \eqref{eq:darcy_stokes_precond} while ($\times$) markers represent
    the (integer-rounded) mean iteration count.
  }
  \label{fig:darcy_stokes_pcg}
\end{figure}

In \Cref{fig:darcy_stokes_2d} we observe that the number of FGMRes iterations is bounded in $\mu$, $K$ and mesh size. Moreover,
using RA for approximating the inverse of $S_{\Gamma}=KL^{1/2}+\mu^{-1}L^{-1/2}$ leads to practically the same convergence as when
the (exact) spectral realization is used. In \Cref{fig:darcy_stokes_pcg} we finally consider behaviour of PCG iterations in
applications of the the Darcy-Stokes preconditioner. When \eqref{eq:dd_preconditioner} is computed using LU
for $(A^{00}_{\Omega})^{-1}$ and \eqref{eq:gevp} in inverting $S_{\Gamma}$, convergence of the
Krylov solver is stable with respect to material and discretization parameters. In case (i) a slight increase of iterations
for large values of $\mu$ can be observed. We recall that for $(A^{00}_{\Omega})^{-1}$ a single AMG V-cycle is applied and
$\epsilon_{\text{RA}}=10^{-12}$ in approximating the Schur component inverse. This observation suggests that sufficient
accuracy is needed when approximating the blocks in \eqref{eq:dd_preconditioner} in order to retain mesh-independent
convergence. This effect will be investigated in future work.

\begin{acknowledgement}
  The author is grateful to prof. Ludmil T. Zikatanov (Penn State) and prof. Kent-Andr{\'e} Mardal
  (University of Oslo) for stimulating discussions on non-overlapping domain decomposition which
  inspired the presented approach. This work received support from the Norwegian Research Council
  grant 303362.
\end{acknowledgement}

\bibliographystyle{spmpsci}
\bibliography{references}

\end{document}